\documentclass{amsart}
\usepackage{amssymb,amscd,amsthm,verbatim}
\input epsf                 

\newtheorem{proposition}{Proposition}[section]
\newtheorem{theorem}[proposition]{Theorem}

\newtheorem{lemma}[proposition]{Lemma}

\newtheorem{prop}[proposition]{Proposition}
\newtheorem{thm}[proposition]{Theorem}
\newtheorem{cor}[proposition]{Corollary}

\newcommand{\naturals}{\mathbb N}

\newcommand{\reals}{\mathbb R}

\newcommand{\Q}{{\mathcal Q}}
\newcommand{\A}{{\mathcal A}}
\newcommand{\B}{{\mathcal B}}
\newcommand{\D}{{\mathcal D}}

\newcommand{\Po}{{\mathcal P}}

\newcommand{\rank}{{\mbox{rank}}}
\newcommand{\0}{{\hat{0}}}

\newcommand{\set}[1]{{\lbrace #1 \rbrace}}

\newcommand{\join}{\vee}

\newcommand{\Span}{\mbox{{\rm Span}}}

\newcommand{\Subcrit}{\mbox{{\rm Subcrit}}}

\newcommand{\nth}{{\mbox{{\small th}}}}

\begin{document}
\title[Dimension of the Poset of Regions]{The Order Dimension of the Poset of Regions in a Hyperplane Arrangement}

\author{Nathan Reading}
\address{
Mathematics Department\\
       University of Michigan\\
       Ann Arbor, MI 48109-1109\\
USA}
\thanks{The author was partially supported by NSF grant DMS-0202430.}
\email{nreading@umich.edu}
\urladdr{http://www.math.lsa.umich.edu/$\sim$nreading/}

\begin{abstract}
We show that the order dimension of the weak order on a Coxeter group of type A, B or D is equal to the rank of the Coxeter group, and give 
bounds on the order dimensions for the other finite types.
This result arises from a unified approach which, in particular, leads to a simpler treatment of the previously known cases, types A and 
B~\cite{Flath,hyperplane}.
The result for weak orders follows from an upper bound on the dimension of the poset of regions of an arbitrary hyperplane arrangement.
In some cases, including the weak orders, the upper bound is the chromatic number of a certain graph.
For the weak orders, this graph has the positive roots as its vertex set, and the edges are related to the pairwise inner products of the roots.
\end{abstract}

\maketitle

\section{Introduction}
\label{main results}

For a finite Coxeter group $W$, let $\dim(W)$ be the order dimension of the weak order on $W$, or in other words the order dimension of the 
poset of regions of the corresponding Coxeter arrangement.
The order dimension of a finite poset $P$ is the smallest $n$ so that $P$ can be embedded as an induced subposet of the componentwise order
on $\reals^n$.

\begin{theorem}
\label{dimensions}
The order dimension of the weak order on an irreducible finite Coxeter group has the following value or bounds:
\[\begin{array}{rcccl}
                    &&    	\dim(A_n)     	& = &	n\\
                    &&	\dim(B_n)     	& = &	n\\
                    &&	\dim(D_n)     	& = &	n\\
              6	&\le &	\dim(E_6)     	& \le &	9\\
              7	&\le & 	\dim(E_7)     	& \le  &	11\\
              8	&\le  & 	\dim(E_8)     	& \le  &	19\\
              4	&\le  &	\dim(F_4)     	& \le  &	5\\
                     &&   	\dim(H_3)     	& = &	3\\
              4	&\le  & 	\dim(H_4)     	& \le  &	6\\
                        &&	\dim(I_2(m))   	& =&	 2
\end{array}\]
\end{theorem}
The order dimension of the weak order on a reducible finite Coxeter group is the sum of the dimensions of the irreducible components. 
The result for $A_n$ was proven previously by Flath~\cite{Flath} using the combinatorial interpretation of $A_n$, while the results for 
$A_n$ and $B_n$ were obtained previously by an argument using supersolvability~\cite{hyperplane}.
Theorem \ref{dimensions} gives values or bounds for all types of finite Coxeter groups, including new results on type D and the exceptional groups.
The theorem is based on a unified approach which, in particular, provides a significantly simpler proof of the results for types A and B.

The lower bounds of Theorem~\ref{dimensions} are easily proven by considering the atoms and coatoms of the posets (Proposition~\ref{lower}).
The upper bounds are proven by way of a more general theorem giving an upper bound on the order dimension of the poset of regions of any 
hyperplane arrangement.
Specifically, for a hyperplane arrangement~$\A$ and a fixed region~$B$, let $\Po(\A,B)$ be the poset of regions, that is, the adjacency graph of the regions of~$\A$,
directed away from~$B$.
Then there is a directed graph $\D(\A,B)$ whose vertex set is~$\A$ such that the following holds:
\begin{theorem}
\label{acyclic}
For a central hyperplane arrangement~$\A$ with base region~$B$, the order dimension of $\Po(\A,B)$ is bounded above by the size of 
any covering of $\D(\A,B)$ by acyclic induced sub-digraphs.
\end{theorem}
By a covering of $\D(\A,B)$ by acyclic induced sub-digraphs we mean a partition $\A=I_1\cup I_2\cup\cdots\cup I_k$ such that each $I_j$ induces 
an acyclic sub-digraph of $\D(\A,B)$.  The size of such a covering is $k$.
It is well-known that, in general, order dimension can be characterized as a problem of covering a directed graph by acyclic induced sub-digraphs (see for example~\cite{hyperplane}).
However, if one does this for $\Po(\A,B)$ one generally gets a directed graph with many more vertices than $\D(\A,B)$.

For a large class of arrangements, the minimal cycles in $\D(\A,B)$ have cardinality two.
Thus the order dimension of $\Po(\A,B)$ is bounded above by the chromatic number of the graph $G(\A,B)$ whose vertex set is $\A$ and 
whose edges are the two-cycles of $\D(\A,B)$.
Other connections between graph coloring and order dimension have been made, for example in~\cite{Fels-Trot,FHRT,Yan}.

When~$\A$ is a Coxeter arrangement, the edges of $G(\A,B)$ can be determined by considering inner products of pairs of roots in the corresponding root system.
This leads to straightforward colorings of the graphs for Coxeter arrangements of types A, B and D.
The dimension results in types G and I are trivial using Theorem~\ref{acyclic} or by much simpler considerations.
The value and bounds for types E, F and H come from computer computations of $\chi(G(\A,B))$.
The programs used for these computations were written by John Stembridge, and are available on the author's website.

The proof of Theorem~\ref{acyclic} uses a new formulation of order dimension, similar in spirit to the formulation in terms of critical pairs~\cite{Rab-Riv}.
A well-known theorem of Dushnik and Miller~\cite{Du-Mil} says that the order dimension of a poset $P$ is the smallest $d$ so that $P$ can be embedded as an 
induced subposet of $\reals^d$.
The components of the embedding need not be linear extensions of $P$, but rather are order-preserving maps of $P$ into linear orders.
Proposition~\ref{box} uses subcritical pairs (see~\cite{hyperplane}) to give conditions on a set of order-preserving maps from $P$ into linear 
orders, which are necessary and sufficient for the maps to be the components of an embedding.
The subcritical pairs of $\Po(\A,B)$ are identified with the shards of $(\A,B)$.
Introduced in~\cite{hyperplane}, the shards are the components of hyperplanes in~$\A$ which result from ``cutting'' the hyperplanes in a certain way.
This geometric information about the subcritical pairs leads to the proof of Theorem~\ref{acyclic}.

Hyperplane arrangements are dual to zonotopes, and the Hasse diagram of $\Po(\A,B)$ is the same as the 1-skeleton of the corresponding zonotope.
Thus, given~$\A$ and~$B$, one might hope to give an embedding of $\Po(\A,B)$ by mapping each region to the corresponding vertex of an 
equivalent zonotope.
We show that this can be done when~$\A$ is a supersolvable arrangement.

The body of the paper is organized as follows.
In Section~\ref{arr} we give definitions and preliminary results about hyperplane arrangements and posets of regions.
Section~\ref{dim} contains background information about order dimension, and states and proves the reformulation mentioned above 
(Proposition~\ref{box}).
Theorem~\ref{acyclic} is proven in Section~\ref{dim po}, while Section~\ref{coxeter} contains the details of the coloring problem in the 
case of Coxeter arrangements, leading to the proof of Theorem~\ref{dimensions}.
Section~\ref{zonotopal} is a discussion of zonotopal embeddings, and Section~\ref{supersolvable} is an application of Sections \ref{dim po} 
and \ref{zonotopal} to the case of supersolvable arrangements.

\section{Hyperplane Arrangements}
\label{arr}
In this section we give definitions related to hyperplane arrangements, and prove some basic facts about join-irreducible and meet-irreducible elements of the 
poset of regions of an arrangement.
An {\em arrangement}~$\A$ is a finite, nonempty collection of {\em hyperplanes} (codimension~1 linear subspaces) in $\reals^n$.
In general, one might consider arrangements of affine hyperplanes, but in this paper all arrangements will consist of hyperplanes containing the 
origin.
Such arrangements are called {\em central}.
The complement of the union of the hyperplanes is disconnected, and the closures of its connected components are called {\em regions}.
The {\em span} of~$\A$, written $\Span(\A)$, is understood to mean the linear span of the normal vectors of~$\A$, and the {\em rank} 
of~$\A$ is the dimension of $\Span(\A)$.

The {\em poset $\Po(\A,B)$ of regions} of~$\A$ with respect to a fixed region~$B$ is a partial order on the regions defined as follows.
Define $S(R_1,R_2)$ to be the set of hyperplanes separating $R_1$ from $R_2$. 
For any region $R$, the set $S(R):=S(R,B)$ is called the {\em separating set} of $R$.
The poset of regions is a partial order on the regions with $R_1\le R_2$ if and only if $S(R_1)\subseteq S(R_2)$.
The fixed region~$B$, called the {\em base region}, is the unique minimal element of $\Po(\A,B)$.
The definition of $\Po(\A,B)$ is an embedding into a product of $|\A|$ chains, so the dimension of $\Po(\A,B)$ is at most $|\A|$.
For more details on this poset, see~\cite{BEZ,Edelman}.

\begin{figure}[ht]
\caption{A hyperplane arrangement~$\A$ with base region~$B$ and the poset of regions $\Po(\A,B)$.}
\label{ex}
\centerline{\epsfbox{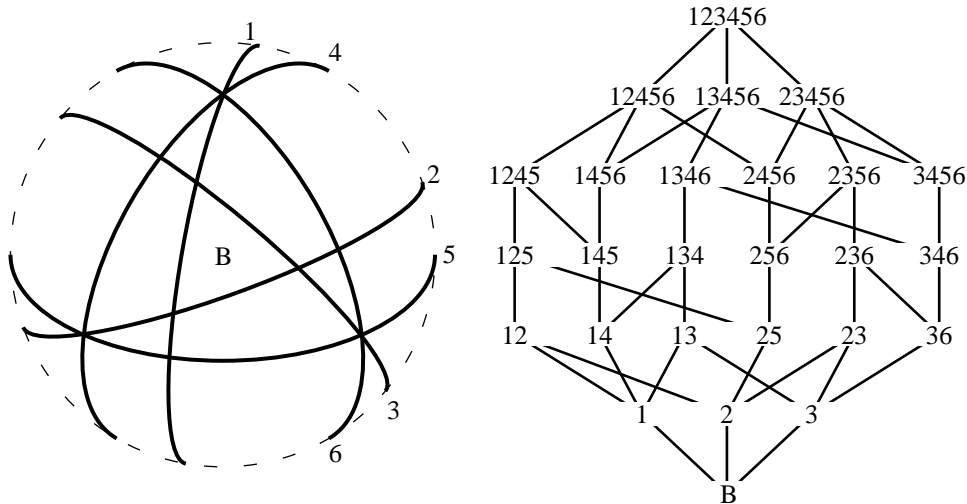}}   
\end{figure}

When~$\A$ is central, the antipodal anti-automorphism of $\Po(\A,B)$, denoted by $R\mapsto -R$, corresponds to complementation of 
separating sets.
In particular there is a unique maximal element $-B$.
A central arrangement is {\em simplicial} if every region is a simplicial cone.
Figure~\ref{ex} shows $\Po(\A,B)$ for a non-simplicial arrangement~$\A$ in $\reals^3$ with base region~$B$.
The hyperplane arrangement is represented as an arrangement of great circles on a 2-sphere.
The northern hemisphere is pictured and the sphere is opaque so that the southern hemisphere is not visible.
The equator is shown as a dotted line to indicate that the equatorial plane is not in~$\A$.
The anti-automorphism $R\mapsto -R$ corresponds to a half-turn of the Hasse diagram of $\Po(\A,B)$.

A subset $\A'\subseteq\A$ is a {\em rank-two subarrangement} if $|\A'|>1$ and there is some codimension-two subspace $L$ of 
$\reals^n$ such that $\A'$ consists of all the hyperplanes containing $L$.
There is a unique region $B'$ of $\A'$ containing~$B$, and the hyperplanes in $\A'$ bounding $B'$ are called {\em basic} hyperplanes in $\A'$.
Rank-two subarrangements and basic hyperplanes are used to define several combinatorial structures which are central to the results in this paper.
The {\em basic digraph} $\D(\A,B)$ is the directed graph whose vertex set is~$\A$, with directed edges $H_1\rightarrow H_2$ whenever $H_1$ 
is basic in the rank-two subarrangement determined by $H_1\cap H_2$.

If $H_1$ and $H_2$ are basic in $\A'$ but $H\in\A'$ is not, then $(H\cap B') = (H_1\cap H_2\cap B')$.
Intersecting both sides of the equality with~$B$, we obtain the following, which we name as a lemma for easy reference later.
\begin{lemma}
\label{basic containment}
If $H_1$ and $H_2$ are basic in $\A'$ but $H\in\A'$ is not, then $(H\cap B) = (H_1\cap H_2\cap B)$.
\qed
\end{lemma}

The bound of Theorem~\ref{acyclic} is not sharp.
For example, an arrangement~$\A$ is {\em 3-generic} if every rank-two subarrangement contains exactly two hyperplanes~\cite{Ziegler}.
For a 3-generic arrangement, $\D(\A,B)$ is complete, in the sense that every pair of vertices is connected by one directed edge in each direction.
Thus Theorem~\ref{acyclic} gives the upper bound $|\A|$ on the order dimension of $\Po(\A,B)$.
There is a unique (up to combinatorial isomorphism) 3-generic arrangement in $\reals^3$ with $|\A|=4$.
The intersection of this arrangement with the unit sphere cuts the sphere into 8 triangles and 6 quadrilaterals, so as to be combinatorially 
isomorphic to the boundary of the cuboctahedron.
If $B$ is chosen to be one of the triangular regions, then $\Po(\A,B)$ has order dimension 3, as can be seen by modifying the usual 
embedding of the Boolean algebra.
In light of Proposition~\ref{lower} which will be proved in Section~\ref{dim}, this example also illustrates the fact that the order dimension depends on the 
choice of base region.

In the example of Figure~\ref{ex}, the rank-two subarrangements are the following subsets of~$\A$:  12, 13, 23, 15, 26, 34, 146, 245 and 356.
Figure~\ref{ex2} shows the basic digraph for this example.
Note the three-cycle $4\rightarrow 5\rightarrow 6\rightarrow 4$.

\begin{figure}[ht]
\caption{The basic digraph $\D(\A,B)$ for $(\A,B)$ as in Figure~\ref{ex}.}
\label{ex2}
\centerline{\epsfbox{hplane_dim_fig2.ps}}   
\end{figure}

The {\em shards} of an arrangement are pieces of the hyperplanes which arise as follows.
For each $H\in\A$, and for each rank-two subarrangement $\A'$ containing $H$, if $H$ is not basic in $\A'$, cut $H$ by removing $L$ from 
$H$, where $L$ is the codimension-two subspace defining $\A'$.
Each hyperplane may be cut several times, and the resulting connected components of the hyperplanes in~$\A$ are called the {\em shards} of~$\A$
with respect to~$B$.
Shards were introduced in~\cite{hyperplane} in connection with certain lattice properties of $\Po(\A,B)$ for a simplicial arrangement~$\A$.

\begin{figure}[ht]
\caption{The decomposition of~$\A$ into shards for $(\A,B)$ as in Figure~\ref{ex}.}
\label{ex3}
\centerline{\epsfbox{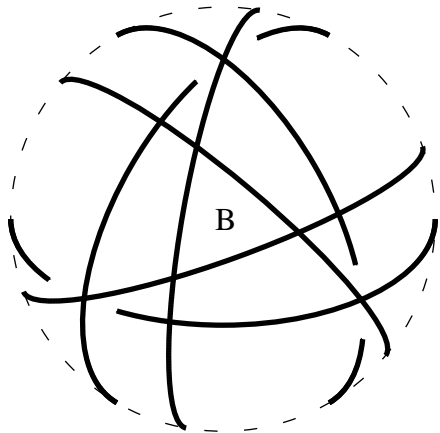}}   
\end{figure}

Figure~\ref{ex3} shows the decomposition into shards of the example $(\A,B)$ of Figures~\ref{ex} and~\ref{ex2}.
Once again, the drawing shows the northern hemisphere.
The southern-hemisphere picture is similar, and in this example all of the shards intersect both hemispheres.

Let $P$ be a poset.
The join $\join X$ of a set $X\subseteq P$ is the unique minimal upper bound for $X$ in $P$, if such exists.   
An element $j$ of a poset $P$ is {\em join-irreducible} if there is no set $X\subseteq P$ with $j\not\in X$ and $j=\join X$.
If $P$ has a unique minimal element $\0$, then $\0$ is $\join\emptyset$ and thus is not join-irreducible.  
Meet-irreducible elements are defined dually.

In a lattice, $j$ is join-irreducible if and only if it covers exactly one element, but this need not be the case in a non-lattice.
However, a region $J$ in $\Po(\A,B)$ is join-irreducible if and only if it covers exactly one region $J_*$, because cover relations in $\Po(\A,B)$ 
correspond to deleting one element from the separating set.
If~$\A$ is a central arrangement, a region $M$ is meet-irreducible if and only if it is covered by exactly one element, denoted $M^*$.
The shards of a finite central arrangement are related to the join- and meet-irreducibles of the poset of regions, as explained below.
Given a shard $\Sigma$, let $H_\Sigma$ be the hyperplane of~$\A$ containing $\Sigma$.
Let $U(\Sigma)$ be the set of {\em upper regions} of $\Sigma$, that is, the set of regions $R$ of~$\A$ which intersect $\Sigma$ in codimension
one and which have $H_\Sigma\in S(R)$.
The set $L(\Sigma)$ of {\em lower regions} of $\Sigma$ is the set of regions $R$ of~$\A$ which intersect $\Sigma$ in codimension one and which 
have $H_\Sigma\not\in S(R)$.
In the following propositions, $U(\Sigma)$ and $L(\Sigma)$ are considered to be subposets of $\Po(\A,B)$.

\begin{prop}
\label{j sigma}
A region $J$ is join-irreducible in $\Po(\A,B)$ if and only if $J$ is minimal in $U(\Sigma^J)$ for some shard $\Sigma^J$, 
in which case $S(J_*)=S(J)-\set{H_{\Sigma^J}}$.
\end{prop}
\begin{proof}
Suppose $J$ is join-irreducible.
Then $J$ and $J_*$ are separated by some shard $\Sigma$ and $S(J_*)=S(J)-\set{H_\Sigma}$.
Since $J$ covers only $J_*$ and $H_\Sigma\not\in S(J_*)$, any region $R<J$ has $H_\Sigma\not\in S(R)$.
In particular, $R$ is not in $U(\Sigma)$, so the region $J$ is minimal in $U(\Sigma)$.
Conversely, suppose $J$ is minimal in $U(\Sigma)$ for some shard $\Sigma$, and suppose that $J$ covers more than one region.
Let $J_*$ be the region whose separating set is $S(J)-\set{H_\Sigma}$.
If $b$ is some vector in $B$, then the facets of $J$ which one would cross to go down by a cover in $\Po(\A,B)$ are the facets of $J$ whose
outward-directed normals have positive inner product with $b$.
In particular, this set of facets is a ball, and therefore we can find a region $R$ covered by $J$ so that $R\cap J\cap J_*$ has codimension two.
Let $S(J)-S(R)=\set{H}$ and let $\A'$ be the rank-two subarrangement containing $H$ and $H_\Sigma$. 
The subarrangement $\A'$ and the regions adjacent to $\cap\A'$ are depicted in Figure~\ref{prooffig}.
\begin{figure}[!h]
\caption{}
\label{prooffig}
\centerline{\epsfbox{hplane_dim_proof_fig.ps}}   
\end{figure}
Since $J$ covers both $J_*$ and $R$ by respectively crossing $H_\Sigma$ and $H$, the hyperplanes $H_\Sigma$ and $H$ are basic in $\A'$.
Because $J$ intersects $\cap\A'$ in codimension two, there is a region $R'$ whose separating set is $(S(J)-\A')\cup H_\Sigma$.
This region is in $U(\Sigma)$, contradicting the minimality of~$J$.
\end{proof}

The following proposition is dual to Proposition~\ref{j sigma}.
\begin{prop}
\label{m sigma}
A region $M$ is meet-irreducible in $\Po(\A,B)$ if and only if $M$ is maximal in $L(\Sigma_M)$ for some shard $\Sigma_M$, 
in which case $S(M^*)=S(M)\cup\set{H_{\Sigma_M}}$.
\qed
\end{prop}
We will write $H^J$ for $H_{\Sigma^J}$ and $H_M$ for $H_{\Sigma_M}$.

We conclude the section with a technical observation which is used in the proof of Theorem~\ref{acyclic}.

\begin{lemma}
\label{shards inherit}
Let~$\A$ be a central hyperplane arrangement with base region~$B$ and let $I\subseteq \A$.
Let $H\in I$ be a sink in the sub-digraph of $\D(\A,B)$ induced by $I$, let $\A^-:=\A-\set{H}$ and let $B^-$ be the region of $\A^-$ containing~$B$.
Then the shards of $(\A,B)$ contained in hyperplanes in $I-\set{H}$ are exactly the shards of $(\A^-,B^-)$ contained in hyperplanes $I-\set{H}$.
\end{lemma}
\begin{proof}
Since $H$ is a sink in the sub-digraph of $\D(\A,B)$ induced by $I$, for any $H'\in I-\set{H}$, the hyperplane $H$ is not basic in the rank-two subarrangement determined by $H\cap H'$.
In particular, removing $H$ has no effect on the process of ``cutting'' $H'$ into shards.
\end{proof}

\section{Order dimension and subcritical pairs}
\label{dim}
In this section we give background information on order dimension and a new formulation of order dimension in terms of subcritical pairs.

A poset $E$ on the same ground set as $P$ is called an {\em extension} of $P$ if $a\le_Pb$ implies $a\le_Eb$.
An extension is called {\em linear} if it is a total order.
The {\em order dimension} $\dim(P)$ of a finite poset $P$ is the smallest $d$ so that $P$ can be written as the intersection---as relations---of 
$d$ linear extensions of $P$.
Say $Q$ is a(n) {\em (induced) subposet} of $P$ if there is a one-to-one map $i:Q\rightarrow P$ such that $x\le_Qy$ if and only if $i(x)\le_Pi(y)$.
If $Q$ is an induced subposet of $P$, then $\dim(Q)\le\dim(P)$.

The ``standard example'' of a poset of dimension $n$ is the collection of subsets of $[n]$ having cardinality 1 or $n-1$.
In an arbitrary finite central arrangement~$\A$ with base region~$B$, the collection of regions covering~$B$ or covered by $-B$ form a subposet of $\Po(\A,B)$
which is isomorphic to a standard example.
Each facet (maximal face) of~$B$ corresponds to a region covering~$B$, and thus we have the following lower bound on $\dim(\Po(\A,B))$.
\begin{prop}
\label{lower}
The order dimension of $\Po(\A,B)$ is at least the number of facets of~$B$, which is at least the rank of $\A$.
\qed
\end{prop}

A pair $(j,m)$ in a poset $P$ is called {\em subcritical} if:
\begin{enumerate}
\item[(i) ] $j\not\le m$,
\item[(ii) ] For all $x\in P$, if $x<j$ then $x\le m$,
\item[(iii) ] For all $x\in P$, if $x>m$ then $x\ge j$.
\end{enumerate}
The set of subcritical pairs of $P$ is denoted $\Subcrit(P)$.
The more commonly used {\em critical pairs} are defined by replacing condition (i) with
\begin{enumerate}
\item[(i') ] $j$ is incomparable to $m$.
\end{enumerate}
Thus critical pairs are in particular subcritical, and a subcritical pair $(j,m)$ that is not critical has the property that $j$ covers $m$ but covers nothing else, and 
$m$ is covered by $j$ and by nothing else.

The following proposition was proven in~\cite{Rab-Riv} for critical pairs in a lattice, and the proof for subcritical pairs in a poset is essentially the same.
\begin{prop}
\label{irr}
If $(j,m)$ is a subcritical pair in a poset $P$, then $j$ is join-irreducible and $m$ is meet-irreducible.
\qed
\end{prop}

An extension $E$ of a poset $P$ is said to {\em reverse} a critical or subcritical pair $(j,m)$ if $m<j$ in $E$.
The following formulation of order-dimension is due to Rabinovitch and Rival.
\begin{prop}\cite{Rab-Riv}
\label{critical}
The order dimension of a finite poset $P$ is equal to the smallest $d$ such that there
exist linear extensions $L_1,\ldots,L_d$ such that for each critical pair $(j,m)$ of $P$ there is some $L_i$ which reverses $(j,m)$.
\qed
\end{prop}
Since critical pairs are in particular subcritical, one can substitute ``subcritical'' for ``critical'' in Proposition~\ref{critical}. 
Subcritical pairs also occur in~\cite{hyperplane}.

A well-known theorem of Dushnik and Miller~\cite{Du-Mil} says that the order dimension of a poset $P$ is the smallest $d$ so that $P$ can be embedded as an 
induced subposet of $\reals^d$.
For a poset $P$ with $|P|=n$ and $\dim(P)=d$, we can use $d$ linear extensions whose intersection is $P$ to embed $P$ as a subposet of $[n]^d$.
The theorem of Dushnik and Miller suggests that we can embed $P$ into a smaller $d$-dimensional ``box.''
Subcritical pairs are the key to embedding a poset into a small box.
Let $\eta:P\rightarrow Q$ be an order-preserving map from $P$ to $Q$.
That is, whenever $x\le y$ in $P$, then $\eta(x)\le\eta(y)$ in $Q$.
Say $\eta$ {\em reverses} a subcritical pair $(j,m)$ if $\eta(m)<\eta(j)$.
The strict inequality is essential here.

\begin{prop}
\label{box}
The order dimension of a finite poset $P$ is equal to the smallest $d$ such that there exist order-preserving maps 
$\eta_1,\ldots,\eta_d:P\to\naturals$ such that for each subcritical pair $(j,m)$ of $P$ there is some $\eta_i$ which reverses $(j,m)$.
\end{prop}
\begin{proof}
Suppose ${\mathbf \eta}=(\eta_1,\eta_2,\ldots,\eta_d)$ is an embedding of $P$ into $\naturals^d$, and let $(j,m)$ be a subcritical pair.
Since $j\not\le m$, there must be some $\eta_i$ which reverses $(j,m)$.
Conversely, suppose that there exist order-preserving maps $\eta_1,\ldots,\eta_d:P\to\naturals$ such that for each subcritical pair $(j,m)$ of $P$ 
there is some $\eta_i$ which reverses $(j,m)$.
Let ${\mathbf \eta}:=(\eta_1,\eta_2,\ldots,\eta_d)$.
To show that ${\mathbf \eta}$ is an embedding, we must show that for any pair $(a,b)$ in $P$ with $a\not\le b$, there is some $i\in[D]$ such that 
$\eta_i(b)<\eta_i(a)$.
The simple proof of this fact follows the proof of Proposition~\ref{critical}.
Suppose $(a,b)$ is an exception, or in other words, $\eta_i(b)\ge\eta_i(a)$ for all $i\in D$.
If there exists $a'< a$ such that $a'\not\le b$, replace $a$ by $a'$ to obtain a new pair $(a',b)$, which is also an exception.
(If $\eta_i(b)<\eta_i(a')$ for some $i$, then because $\eta_i$ is order-preserving we have $\eta(b)<\eta(a')\le\eta(a)$, 
contradicting the fact that $(a,b)$ was an exception.)
Similarly, if there exists $b'>b$ with $a\not\le b'$, the pair $(a,b')$ is an exception.
Continue making these replacements, and since $a$ always moves down in the poset and $b$ always moves up, the process will eventually 
terminate by finding an exception which is also a subcritical pair.
This contradiction shows that ${\mathbf \eta}$ is indeed an embedding.
\end{proof}

Some modifications of Proposition~\ref{box} are worth mentioning, although they will not be used in this paper.
Similar modifications of Proposition~\ref{critical} are given in~\cite[Section 1.12]{Trotter}.

\begin{prop}
\label{flubby box}
The order dimension of a finite poset $P$ is the smallest $d$ such that there exist posets $Q_i$ and order-preserving maps $\eta_i:P\rightarrow Q_i$ for $i\in[d]$,
such that for each subcritical pair $(j,m)$ of $P$ there is some $\eta_i$ which reverses $(j,m)$.
\qed
\end{prop}

\begin{prop}
\label{flubbier box}
The order dimension of a finite poset $P$ is the smallest $d$ such that there exist posets $Q_i$, subposets $P_i$ of $P$ and order-preserving maps 
$\eta_i:P_i\rightarrow Q_i$ for $i\in[d]$, such that for each subcritical pair $(j,m)$ of $P$ there is some $i$ with $j,m\in P_i$ and $(j,m)$ reversed 
by $\eta_i$.
\qed
\end{prop}

Proposition \ref{flubby box} follows from Proposition \ref{box} by considering linear extensions of the $Q_i$.
Proposition \ref{flubbier box} follows from Proposition \ref{flubby box} via the following observation:
If $P'$ is an induced subposet of a finite poset $P$, then any order-preserving map $\eta':P'\rightarrow Q$ can be extended to an order preserving 
map $\eta:P\rightarrow E$, where $E$ is some extension of $Q$.

\section{Order dimension of the poset of regions}
\label{dim po}
In this section we relate the shards of $(\A,B)$ to the subcritical pairs in $\Po(\A,B)$.
This relationship, along with Proposition~\ref{box}, is then used to prove Theorem~\ref{acyclic} via an explicit embedding.

\begin{prop}
\label{sc sigma}
Let~$\A$ be a central arrangement.
A pair $(J,M)$ in $\Po(\A,B)$ is subcritical if and only if there is a shard $\Sigma$ such that $J$ is minimal in $U(\Sigma)$, $M$ is maximal in $L(\Sigma)$ and
$J_*\le M$.
\end{prop}
\begin{proof}
Suppose $(J,M)$ is subcritical.
Then by Proposition~\ref{irr}, $J$ is join-irreducible and $M$ is meet-irreducible, so $J_*$ and $M^*$ are defined.
By condition (ii), $J_*\le M$, and in light of Propositions~\ref{j sigma} and~\ref{m sigma} it remains to show that $\Sigma^J=\Sigma_M$.
By condition (iii), $J\le M^*$ as well.
Thus we have $S(J_*)\subseteq S(M)$ and $S(J)\subseteq S(M^*)$.
Therefore $H^J\in S(M^*)$.
If we also have $H^J\in S(M)$ then $S(J)\subseteq S(M)$, contradicting the fact that $(J,M)$ is a subcritical pair.
So $H^J=H_M$, or in other words $\Sigma^J$ and $\Sigma_M$ are contained in the same hyperplane.
Suppose for the sake of contradiction that $\Sigma^J\neq\Sigma_M$.
Then there is a codimension-two subspace $L$ between $\Sigma^J$ and $\Sigma_M$ in $H^J$ such that $H^J$ is not basic in the 
associated rank-two subarrangement $\A'$.
Then necessarily, one of the two basic hyperplanes is in $S(J_*)\cap\A'$ but not in $S(M)\cap\A'$.
This contradiction to $J_*\le M$ shows that $\Sigma^J=\Sigma_M$.

Conversely, suppose that there is a shard $\Sigma$ such that $J$ is minimal in $U(\Sigma)$, $M$ is maximal in $L(\Sigma)$ and $J_*\le M$.
Then by Propositions~\ref{j sigma} and~\ref{m sigma}, $J$ is join-irreducible and $M$ is meet-irreducible and because $J_*\le M$ we have $J\le M^*$ as well.
Thus conditions (ii) and (iii) are satisfied.
Since $M\in L(\Sigma)$ and $J\in U(\Sigma)$, we have $J\not\le M$.
\end{proof}

\begin{lemma}
\label{shard arrow}
Let $(J,M)$ be a subcritical pair in $\Po(\A,B)$ for a central arrangement~$\A$ and let $H\in\A$.
If $H\not\in S(J)$ and $H\in S(M)$, then $H^J$ is basic in the rank-two subarrangement determined by $H\cap H^J$.
\end{lemma}
\begin{proof}
Suppose $H\not\in S(J)$ and $H\in S(M)$ for some critical pair $(J,M)$ and let $\A'$  be the rank-two subarrangement determined by $H\cap H^J$.
By Proposition~\ref{sc sigma}, the codimension-one faces $J\cap J_*$ and $M\cap M^*$ are in the same shard $\Sigma_J$, and thus
in particular $H^J$ is basic in $\A'$.
\end{proof}

Let $I\subseteq\A$ induce an acyclic sub-digraph on $\D(\A,B)$.
Let $F_I$ be the set of subcritical pairs $(J,M)$ in $\Po(\A,B)$ such that $H^J\in I$.
Let $H_1,H_2,\ldots,H_{|I|}$ be an ordering of the hyperplanes in $I$ such that whenever $H_i\rightarrow H_j$ in $\D(\A,B)$, we have $i<j$.
For any region $R$ of~$\A$, let $\eta_I(R)$ be the word of length $|I|$ in 0's and 1's whose $i^\nth$ letter is 0 if $H_i\not\in S(R)$ 
and 1 if $H_i\in S(R)$.
Thinking of this word as a binary number, we have constructed a map $\eta_I$ from $\Po(\A,B)$ to the interval $[0,2^{|I|}-1]$.
The map is order-preserving because the order on $\Po(\A,B)$ is containment of separating sets.

\begin{lemma}
\label{binary}
The map $\eta_I$ reverses all of the subcritical pairs in $F_I$.
\end{lemma}
\begin{proof}
The proof is by induction on $k:=|I|$.
If $k=1$, the result is trivial, so suppose $k>1$, consider the arrangement $\A^-:=\A-\set{H_k}$, with base region $B^-$ as in Lemma 
\ref{shards inherit}.
The hyperplane $H_k$ is a sink in the sub-digraph of $\D(\A,B)$ induced by $I$, so $I^-:=I-\set{H_k}$ induces an acyclic sub-digraph of 
$\D(\A,B)$.
By Lemma~\ref{shards inherit}, the shards of $(\A,\B)$ contained in hyperplanes of $I-\set{H}$ are exactly the shards of $(\A^-,B^-)$ 
contained in hyperplanes of $I-\set{H}$.
The notation $\eta_{I^-}$ could be interpreted either as a map on $\Po(\A,B)$ or on $\Po(\A^-,B^-)$.
However, for a region $R$ of~$\A$, if $R^-$ is the region of $\A^-$ containing $R$, then $S(R,B)\cap I^-=S(R^-,B^-)\cap I^-$, so the 
distinction is meaningless.
If $(J,M)$ is a subcritical pair in $F_I$ not associated with the hyperplane $H_k$, then by Lemma~\ref{shards inherit} and 
Proposition~\ref{sc sigma}, $(J^-,M^-)$ is a subcritical pair in $\Po(\A^-,B^-)$ associated to some hyperplane in $I^-$.
Thus by induction, $\eta_{I^-}(M^-)<\eta_{I^-}(J^-)$.
This is a strict inequality in the lexicographic order, and since $\eta_I$ is obtained from $\eta_{I^-}$ by appending an additional digit on the 
right, the strict inequality is preserved regardless of what the new digits are.
Thus we have $\eta_I(M)<\eta_I(J)$.

If, on the other hand, $(J,M)$ is a subcritical pair associated with $H_k$, the last digit of $\eta_I(J)$ is 1 and the last digit of $\eta_I(M)$ is 0.
Thus if we can show that $S(M)\cap I\subseteq S(J)\cap I$, we will have $\eta_I(M)<\eta_I(J)$.
Suppose to the contrary that there is some $H\in I$ with $H\in S(M)$ but $H\not\in S(J)$.
Then Lemma~\ref{shard arrow} says that $H_k$ is basic in the rank two subarrangement $\A'$ determined by $H$ and $H_k$.
However, this means that $H_k\rightarrow H$ in $\D(\A,B)$, and thus $H$ should have occurred after $H_k$ in the ordering on $I$.
\end{proof}

Recall that Theorem~\ref{acyclic} states that the order dimension of $\Po(\A,B)$ is bounded above by the smallest $k$ such that 
$\A=I_1\cup I_2\cup\cdots\cup I_k$ and $I_j$ induces an acyclic sub-digraph of $\D(\A,B)$ for each $j$.
\begin{proof}[Proof of Theorem~\ref{acyclic}]
Lemma~\ref{binary} can be used for each sub-digraph to obtain the components of an order-preserving map 
$\Po(\A,B)\rightarrow\naturals^k$ which satisfies the hypotheses of Proposition~\ref{box}.
\end{proof}

\begin{figure}[ht]
\caption{An embedding of the poset of regions $\Po(\A,B)$ into $[0,3]^3$ for $(\A,B)$ as in Figure~\ref{ex}.}
\label{ex4}
\centerline{\epsfbox{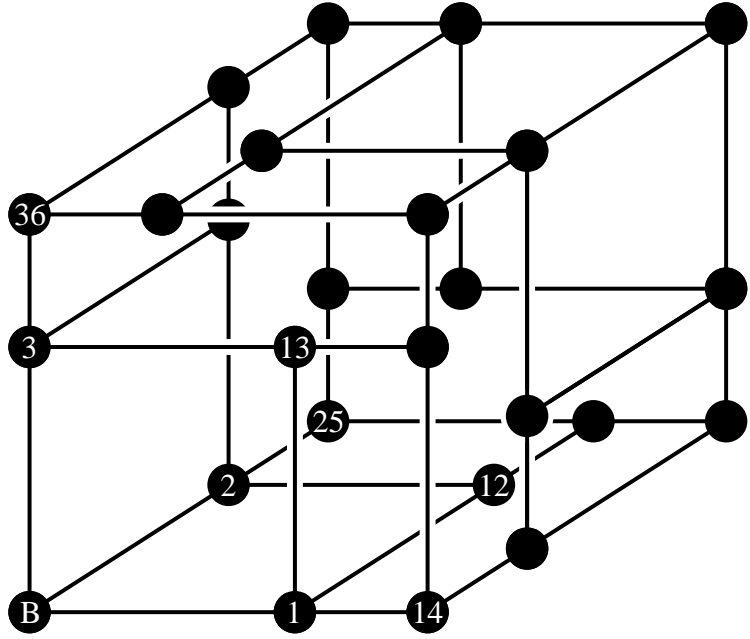}}   
\end{figure}

The directed graph in Figure~\ref{ex2} can be partitioned into three acyclic sub-digraphs, but not fewer.
The partition is $I_1:=\set{1\rightarrow 4}$, $I_2:=\set{2\rightarrow 5}$, $I_3:=\set{3\rightarrow 6}$.
Let $\eta_1:=\eta_{I_1}$ as in Lemma~\ref{binary}, and similarly $\eta_2$ and $\eta_3$.
The image of the map ${\mathbf \eta}=(\eta_1,\eta_2,\eta_3)$ is illustrated in Figure~\ref{ex4}.
In this figure, the first coordinate of ${\mathbf \eta}$ is the horizontal axis, the third coordinate is the vertical axis, and the positive direction of 
the 2nd coordinate points down into the page.
It may also aid the reader's visualization to know that in this example, all of the regions of~$\A$ map to the boundary of the cube.

The {\em basic graph} $G(\A,B)$ is the graph whose vertex set is~$\A$, with edges $\set{H_1,H_2}$ whenever $H_1$ and $H_2$ are the 
basic hyperplanes in some rank-two subarrangement. 
The directed graph $\Q(\A,B)$ has vertex-set~$\A$, with $H\rightarrow H'$ whenever $H$ is a basic hyperplane in some rank-two
subarrangement and $H'$ is a non-basic hyperplane in the same subarrangement.
The edges in $G(\A,B)$ are exactly the directed two-cycles in $\D(\A,B)$.
The directed graph $\Q(\A,B)$ is obtained from $\D(\A,B)$ by deleting the directed edges which are contained in two-cycles.
The following is an immediate corollary of Theorem~\ref{acyclic}.
\begin{cor}
\label{graph}
If $\Q(\A,B)$ is acyclic, then $\dim(\Po(\A,B))\le\chi(G(\A,B))$.
\end{cor}
Here $\chi(G)$ is the chromatic number of the graph $G$.
The acyclicity of $\Q(\A,B)$ also has important consequences for order-theoretic and lattice-theoretic properties of $\Po(\A,B)$ 
\cite{hyperplane}.
In the example of Figures~\ref{ex} through~\ref{ex4}, $\Q(\A,B)$ is not acyclic, and thus Corollary~\ref{graph} does not apply.

\section{Colorings of root systems}
\label{coxeter}
In this section, we use Corollary~\ref{graph} to relate the dimension of the weak order on a finite Coxeter group to a coloring problem on the 
corresponding root system.
Colorings are given which prove Theorem~\ref{dimensions} for types A, B, D and I.  
For types E, F, and H, the bounds were determined using computer programs written by John Stembridge and available on the author's website.

Given a non-zero vector $v$ in $\reals^n$, let $H_v$ be the hyperplane normal to $v$, and let $r_v$ be the Euclidean reflection fixing $H_v$.
A {\em (finite) root system} is a finite collection of vectors in $\reals^n$, satisfying the following properties:
\begin{enumerate}
\item[(i) ] For any $\beta\in\Phi$, we have $r_\beta\Phi=\Phi$. 
\item[(ii) ] For any $\beta\in\Phi$, we have $\beta\reals\cap\Phi=\set{\pm\beta}$.
\end{enumerate}
The group $W$ generated by the reflections $r_\beta$ for $\beta\in\Phi$ is a finite {\em Coxeter group}, and the arrangement of hyperplanes 
$\A_\Phi:=\set{H_\beta:\beta\in\Phi}$ is a {\em Coxeter arrangement}.
Each hyperplane corresponds to two roots.
The rank of a root system $\Phi$ is the dimension of its linear span or equivalently, it is the rank of $\A_\Phi$.
Coxeter arrangements are simplicial, and $W$ acts transitively on the regions of $\A_\Phi$.
Choose some base region~$B$, and for each hyperplane $H$ in $\A_\Phi$, choose the normal root $\beta^+_H$ so that for each region $R$, 
the separating set $S(R)$ is exactly the set of hyperplanes $H$ with $\langle x,\beta^+_H\rangle>0$ for every $x$ in the interior of $R$.
The set $\Phi^+:=\set{\beta^+_H:H\in\A_\Phi}$ is the set of {\em positive roots} of $\Phi$.
Sometimes it is convenient to blur the distinction between the set of hyperplanes and the set of positive roots.
So, for example, we will talk about rank-two subarrangements of root systems, and basic roots in a rank-two subarrangement.

Consider the set of hyperplanes defining facets of~$B$, and call the corresponding set of positive roots the {\em simple roots} $\Delta$.
Since $\A_\Phi$ is simplicial, $\Delta$ is a set of linearly independent vectors.
The set $\set{r_\alpha:\alpha\in\Delta}$ is a set of {\em simple reflections} which generate $W$.
For more details on root systems and Coxeter groups, the reader is referred to~\cite{Bourbaki,Humphreys}.

Root systems have been classified, and we will name Coxeter arrangements according to their corresponding root systems.
There are infinite families $A_n$, $B_n$, $C_n$ and $D_n$, and exceptional root systems $E_6$, $E_7$, $E_8$, $F_4$, $G_2$, $H_3$, $H_4$ and $I_2(m)$.
The root systems $B_n$ and $C_n$ correspond to the same Coxeter arrangement, so we will only consider $B_n$.
Since $G_2$ is the same as $I_2(6)$, we will not consider it separately.
In what follows, we will present specific examples of each type of root system by specifying a set of positive roots.
That set of positive roots determines the associated Coxeter arrangement~$\A$ and the choice of base region~$B$, and for convenience we 
will substitute the name of the root system for the notation $(\A,B)$.
For example, we will refer to $\Po(A_n)$, $\D(A_n)$ and $G(A_n)$ with the obvious meanings.

The poset $\Po(\A_\Phi,B)$ is isomorphic to the {\em weak order} on $W$.
We wish to use root systems to apply Corollary~\ref{graph} to posets of regions of Coxeter arrangements, or equivalently, to the weak orders on the 
corresponding Coxeter groups.
Caspard, Le Conte de Poly-Barbut and Morvan showed that $\Q(\A,B)$ is acyclic whenever~$\A$ is a Coxeter arrangement~\cite{boundedref}.
This was done, using different notation, in the course of establishing a lattice-theoretic result about the weak order on a finite Coxeter group.
Theorem 28 of~\cite{hyperplane} is a different, more geometric proof of the acyclicity of $\Q(\A,B)$ in the case of a Coxeter arrangement.
The acyclicity of $\Q(\A,B)$ allows us to use the more straightforward bound of Corollary~\ref{graph}.
The key to applying Corollary~\ref{graph} is to relate $G(\A,B)$ to the inner products of roots.

If the roots in $\Phi$ consist of more than one $W$ orbit, one can rescale the roots without altering properties (i) and (ii) as long as the rescaling 
is uniform on each $W$-orbit.
For a suitable scaling, the root system has the property that in any rank-two subarrangement, the basic roots are the unique pair of distinct roots 
which minimize the pairwise inner products of distinct positive roots in that rank-two subarrangement.
All of the root systems presented here are scaled so as to have that property.

\subsection*{Type A}
The Coxeter arrangement $A_{n-1}$ corresponds to the root system whose positive roots are 
$\set{\epsilon_i-\epsilon_j:1\le j<i\le n}$.
This root system has rank $n-1$. 
Rank-two subarrangements of the root system $A_{n-1}$ come in two different forms:
A pair of positive roots whose inner product is zero, or a set of three positive roots whose pairwise inner products are 1, 1 and -1.
The hyperplanes corresponding to a pair of orthogonal roots are joined by an edge in $G(A_{n-1})$, and the basic roots in a rank-two 
subarrangement of cardinality three are the pair whose inner product is -1.
Thus independent sets in $G(A_{n-1})$ are sets of roots in which all pairwise inner products are 1.
It is easy to identify the maximal independent sets as having the form $J_i:=\set{\epsilon_i-\epsilon_j:1\le j<i}$ for some fixed $i$ or the form
$\set{\epsilon_i-\epsilon_j:j<i\le n}$ for some fixed $j$.
One $(n-1)$-coloring of $G(A_{n-1})$ uses the sets $J_i$ for $i=2,3,\ldots,n$.

It is also easy to specify the basic digraph $\D(A_{n-1})$.
Besides the 2-cycles, the directed edges are of the form $\epsilon_i-\epsilon_j\rightarrow\epsilon_i-\epsilon_k$  and 
$\epsilon_j-\epsilon_k\rightarrow\epsilon_i-\epsilon_k$ whenever $k<j<i$.
This is because $\set{\epsilon_i-\epsilon_j,\epsilon_i-\epsilon_k,\epsilon_j-\epsilon_k}$ is a rank-two subarrangement whose basic roots are
$\epsilon_i-\epsilon_j$ and $\epsilon_j-\epsilon_k$.

The regions defined by $A_{n-1}$ are in bijection with permutations $\pi_1\pi_2\cdots\pi_n$ of $[n]$.
This notation means that $\pi:i\mapsto \pi_i$.
The separating set of a region corresponds to the {\em inversion set} $I(\pi):=\set{(i,j):1\le i<j\le n:\pi_i>\pi_j}$, and containment of inversion sets is 
called the {\em weak order} on the symmetric group $S_n$.
Thus the coloring of $G(A_{n-1})$ described above and the maps defined in Lemma~\ref{binary} give  an embedding of the 
weak order on $S_n$ into $\reals^{n-1}$.
Specifically, for $i=2,3,\ldots,n$, let $\eta_i(\pi):=\set{j:j<i,\pi_j>\pi_i}$, and interpret this set as a binary number by letting $j$ correspond to the 
$j^\nth$ digit.
This is an embedding of the weak order on $S_n$ into the product $[0,1]\times[0,3]\times[0,7]\times\cdots\times[0,2^{n-1}-1]$.

\subsection*{Type B}
The root system $B_n$ has positive roots \[\set{\epsilon_i\pm\epsilon_j:1\le j<i\le n}\cup\set{\epsilon_i:i\in[n]}.\]
Rank-two subarrangements of $B_n$ can consist of two or three positive roots with the same pairwise inner products as in type A,
or they can be a set of four positive roots whose pairwise inner products are -1, 0, 0, 1, 1, and 1.
The edges in $G(B_n)$ are pairs of roots with inner product -1 and some pairs of roots which have inner product zero.
The rank-two subarrangements of cardinality four have the form $\set{\epsilon_i,\epsilon_j,\epsilon_i\pm\epsilon_j}$, with basic roots 
$\epsilon_j$ and $\epsilon_i-\epsilon_j$.
Thus pairs of the form $\set{\epsilon_i,\epsilon_j}$ and $\set{\epsilon_i+\epsilon_j,\epsilon_i-\epsilon_j}$ are non-edges in $G(B_n)$ even though
these pairs have inner product zero.

Noting that $B_1\subseteq B_2\subseteq\cdots\subseteq B_n$, we obtain an $n$-coloring by setting $I_j=B_j-B_{j-1}$ for $j=1,2,\ldots,n$. 
Another particularly nice coloring decomposes the positive roots into colors of size $n$ so that any pair of roots in the same color have inner 
product 1.
The $i^\nth$ color in this coloring is the set
\[\set{\epsilon_i-\epsilon_j:1\le j<i}\cup\set{\epsilon_i}\cup\set{\epsilon_i+\epsilon_k:1\le i<k}.\]
Using these two colorings, one constructs maps, as in Lemma~\ref{binary}, to embed $\Po(B_n)$ into 
$[0,1]\times[0,7]\times\cdots\times[0,2^{2n-1}-1]$ or into $[0,2^n-1]^n$

\begin{figure}[ht]
\caption{Two colorings of $G(B_3)$.}
\label{exB}
\centerline{\epsfbox{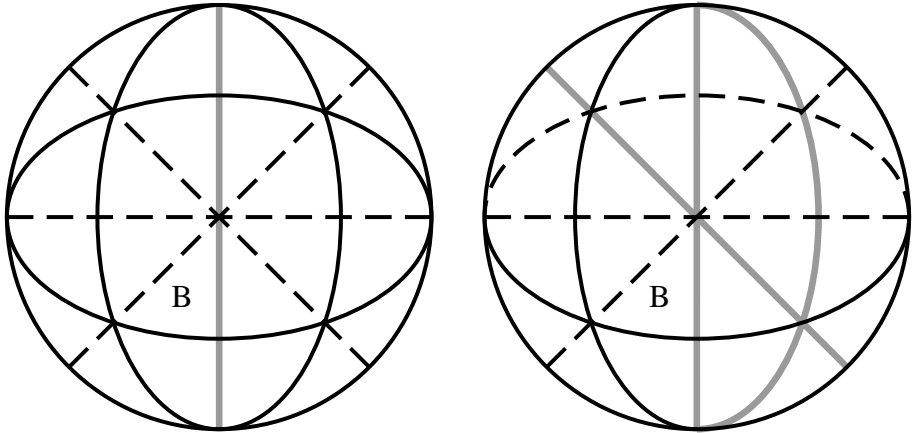}}   
\end{figure}

Figure~\ref{exB} shows these two colorings of $G(B_3)$, the basic graph of the Coxeter arrangement $B_3$.
In this figure, the vector $\epsilon_1$ points to the right, $\epsilon_2$ points towards the top of the page, and $\epsilon_3$ points down into the page.
The hyperplanes are colored in three colors: black, gray and dotted.

\subsection*{Type D}
The positive roots of $D_n$ are $\set{\epsilon_i\pm\epsilon_j:1\le j<i\le n}$.
Rank-two subarrangements of the $D_n$ consist of two or three positive roots with the same pairwise inner products as in type A, so 
the edges in $G(D_n)$ are pairs of roots with inner product -1.
One can color the positive roots by restricting the second coloring given above for $B_n$.
Specifically, the $i^\nth$ color is the set
\[\set{\epsilon_i-\epsilon_j:1\le j<i}\cup\set{\epsilon_i+\epsilon_k:1\le i<k}.\]
This gives an embedding of $\Po(D_n)$ into $[0,2^{n-1}-1]^n$.

\subsection*{Type I}
The graph $G(I_2(m))$ has only a single edge, and thus is two-colorable.
It is also readily apparent by inspection that the dimension of $\Po(I_2(m))$ is two.

\subsection*{Other types}
In each of the infinite families of Coxeter arrangements, the upper bound from Corollary~\ref{graph} agrees with the lower bound of 
Proposition~\ref{lower}, and thus the order dimension equals the rank of the arrangement.
Intriguingly, the situation is different for most of the exceptional groups.
The computational results are:
\begin{eqnarray*}
                        & \chi(G(E_6))     & =  9\\
                        & \chi(G(E_7))     & =  11\\
           16   \le & \chi(G(E_8))     & \le  19\\
                        & \chi(G(F_4))     & =  5\\
                        & \chi(G(H_3))     & =3  \\
                        & \chi(G(H_4))     & = 6\\
\end{eqnarray*}
Of the six Coxeter arrangements of types E, F and H, only $H_3$ has the property that the chromatic number of $G$ is equal to the rank
of the arrangement.

\section{Zonotopal embeddings}
\label{zonotopal}
In this section, we define zonotopal embeddings of the poset of regions, and prove a proposition which gives sufficient conditions for constructing
such embeddings.
In Section \ref{supersolvable}, we apply these condition to supersolvable arrangements.

Given an arrangement~$\A$ and base region~$B$, one can choose a set of normal vectors $\set{n_H:H\in\A}$ such that for each region $R$, 
the separating set $S(R)$ is exactly the set of hyperplanes $H$ with $\langle x,n_H\rangle>0$ for every $x$ in the interior of $R$.
One associates a zonotope to $(\A,B)$ by taking the Minkowski sum of the line segments connecting the origin to each $n_H$.
The 1-skeleton of this zonotope, directed away from the origin, defines a poset isomorphic to $\Po(\A,B)$.
The isomorphism is $Z~:~R\mapsto \sum_{H\in S(R)}n_H$.
The combinatorial type of the zonotope (and thus the partial order) is not changed when the normal vectors are scaled by positive constants.

One might hope that, with some suitable scaling of the normals, and some choice of basis for $\reals^n$, the map $Z$ is an embedding (in the sense of
order-dimension) of $\Po(\A,B)$ into $\reals^n$.
Specifically, choose a basis $b_1,b_2,\ldots,b_n$ for $\reals^n$, and for any vector $v\in\reals^n$, let $v_i$ be the coefficient of $b_i$ when $v$ is expanded
in terms of the basis $b_1,b_2,\ldots,b_n$.
Let $Z_i$ be the map $R\mapsto(Z(R))_i$, the $i^\nth$ component of the vector $Z(R)$.
Call $Z$ a {\em zonotopal embedding} of $\Po(\A,B)$ if for every pair of regions of~$\A$, we have $R_1\le R_2$ if and only if $Z_i(R_1)\le Z_i(R_2)$ for all 
$i\in [n]$.

As an example, consider the hyperplane arrangement in $\reals^2$ whose normal vectors are $n_{H_1}=(1,0)$, $n_{H_2}=(0,1)$ and $n_{H_3}=(1,1)$, choose~$B$ to be the
region containing the vector $(-1,-1)$, and let the $b_i$ be the standard basis.
In this case $Z$ is not an embedding in the sense of order dimension.
Consider the regions $R_1$ and $R_2$ with $S(R_1)=\set{H_1,H_3}$ and $S(R_2)=\set{H_2}$.
We have $R_1\not\ge R_2$, but $Z(R_1)=(2,1)>(0,1)=Z(R_2)$.
However, we can obtain the same arrangement by choosing $n_{H_2}=(r,r)$ for any $r>0$, and when $r<1$, the map $Z$ is a zonotopal embedding.

We now prove a proposition which will help us, in some cases, to find a scaling of the normals so that $Z$ is an embedding.
For each $H\in\A$, define $\nu(H):=\set{H'\in\A:H\rightarrow H' \mbox{ in }\D(\A,B)}$.
Recall that in $\D(\A,B)$, we have $H\rightarrow H'$ whenever $H$ is basic in the rank-two subarrangement $\A'$ determined by $H\cap H'$.
\begin{prop}
\label{sufficient}
Suppose for some $H\in\A$ that
\[(n_H)_i>\sum_{H'\in \nu(H)}(n_{H'})_i.\]
Then the map $Z_i$ reverses all subcritical pairs whose associated hyperplane is $H$.
\end{prop}
\begin{proof}
Let $(J,M)$ be a subcritical pair associated to $H$.
We need to show that $Z_i(M)<Z_i(J)$.
By canceling terms occurring on both sides of the comparison, we see that this is equivalent to proving that
\[\sum_{H'\in S(M)-S(J)}(n_{H'})_i<\sum_{H'\in S(J)-S(M)}(n_{H'})_i.\]
But $H$ is the unique hyperplane in $S(J)-S(M)$, so the right hand sum is $(n_H)_i$.
Any hyperplane $H'$ in $S(M)-S(J)$ intersects the shard associated to $(J,M)$.
If we had $H\not\rightarrow H'$ in $\D(\A,B)$, the intersection $H\cap H'$ would coincide with a cutting of $H$ into shards, and $H'$ would not intersect any 
shard in $H$.
Thus $S(M)-S(J)\subseteq\nu(H)$.
Now we have 
\[(n_H)_i>\sum_{H'\in \nu(H)}(n_{H'})_i\ge\sum_{H'\in S(M)-S(J)}(n_{H'})_i.\]
\end{proof}

\section{Supersolvable arrangements}
\label{supersolvable}
In this section we apply Theorem~\ref{acyclic} and Proposition~\ref{sufficient} to supersolvable arrangements.
The result is a tidier proof of a theorem of~\cite{hyperplane} on the order dimension of the poset of regions of a supersolvable arrangement, and 
a proof that these posets admit zonotopal embeddings.
A Coxeter arrangement is supersolvable if and only if it is of type~A or~B~\cite{Bar-Ihr}, so in particular, weak orders on $A_n$ and $B_n$
admit zonotopal embeddings.

An arrangement~$\A$ is supersolvable if its lattice of intersections $L(\A)$ is supersolvable.
The reader unfamiliar with $L(\A)$ and/or supersolvability can take the following theorem to be the definition of a supersolvable arrangement, 
or see~\cite{BEZ,Or-Ter} for definitions.
\begin{thm}\cite[Theorem~4.3]{BEZ}
\label{characterization}
Every hyperplane arrangement of rank 1 or 2 is supersolvable.
A hyperplane arrangement~$\A$ of rank $d\ge 2$ is supersolvable if and only if it can be written as $\A=\A_0\uplus\A_1$, where
\begin{enumerate}
\item[(i) ]$\A_0$ is a supersolvable arrangement of rank $d-1$.
\item[(ii) ]For any $H',H''\in\A_1$, there is a unique $H\in\A_0$ such that $H'\cap H''\subseteq H$.
\end{enumerate}
\qed
\end{thm}
Here ``$\uplus$'' refers to disjoint union.

Since $\A_0$ has rank one less than~$\A$, the intersection of $\cap\A_0$ with $\Span(\A)$ has dimension 1.
Call this subspace $D$.
\begin{lemma}
\label{no H}
If $H\in\A_1$ then $D\not\subseteq H$.
\end{lemma}
\begin{proof}
Suppose that $D\subseteq H'$ for some $H'\in\A_1$.
Since the rank of~$\A$ is strictly greater than the rank of $\A_0$, there is some $H''\in\A_1$ not containing $D$.
Then $H'\cap H''$ is contained in some unique hyperplane $H$ of $\A_0$.
But then $H=H'$, because both contain the span of $D$ and $H'\cap H''$.
This contradicts the fact that~$\A$ is the disjoint union of $\A_0$ and $\A_1$.
\end{proof}

Let $R$ be a region of $\A_0$, let $v$ be any vector in the interior of $R$.
By Lemma~\ref{no H}, no hyperplane in $\A_1$ contains $D$, so the affine line $v+D$ intersects every hyperplane in $\A_1$.
By Theorem~\ref{characterization}(ii), we can linearly order the hyperplanes of $\A_1$ according to where they intersect $v+D$, and this ordering does not 
depend on the choice of $v\in\mbox{int}(R)$, but only on a choice of direction on $D$.
In particular, consider the set of regions of~$\A$ contained in $R$:  the graph of adjacency on these regions is a path.
As in~\cite{BEZ}, define a {\em canonical base region} inductively:
Any region of an arrangement of rank 2 is a canonical base region.
For a supersolvable arrangement $\A=\A_0\uplus\A_1$, and a region $R$ of~$\A$, let $R_0$ be the region of $\A_0$ containing $R$.
Then~$B$ is a canonical base region if $B_0$ is a canonical base region of $\A_0$ and if the regions of~$\A$ contained in $B_0$ are linearly ordered 
in $\Po(\A,B)$.
The linear order on the regions of~$\A$ contained in $B_0$ also gives a linear order $H_1,H_2,\ldots,H_k$ on the hyperplanes in $\A_1$.

\begin{prop}
\label{superindependent}
If $\A=\A_0\cup\A_1$ is a supersolvable arrangement and~$B$ is a canonical base region, then $\A_1$ induces an acyclic sub-digraph of $\D(\A,B)$.
\end{prop}
\begin{proof}
First we show that there are no 2-cycles in the sub-digraph of $\D(\A,B)$ induced by $\A_1$.
Suppose to the contrary that $H'$ and $H''$ in $\A_1$ are both basic in the rank-two subarrangement $\A'$ they determine.
By Theorem~\ref{characterization}, there is a unique $H\in\A_0\cap\A'$, and Lemma~\ref{basic containment} says that 
$(H\cap B) = (H'\cap H''\cap B)$.
But $H\cap B$ intersects $D$ in dimension one, and thus so does $H'\cap H''\cap B$.
In particular, $H'\cap H''$ contains $D$, contradicting Lemma~\ref{no H}.
This contradiction proves that there are no 2-cycles in the sub-digraph of $\D(\A,B)$ induced by $\A_1$.

Next, we claim that whenever $H_i\rightarrow H_j$ in $\D(\A,B)$, for $H_i,H_j\in \A_1$, we must have $i<j$.
To see this, consider starting at some vector $v$ in the interior of~$B$ and moving along $v+D$ in such a direction as to meet the hyperplanes in $\A_1$.
Since $H_i\rightarrow H_j$, the hyperplane $H_i$ is basic in the rank-two subarrangement $\A'$ determined by $H_i\cap H_j$, and by the previous paragraph, 
no other hyperplane in $\A_1$ is basic in $\A'$.
As we move along $v+D$, we must cross a basic hyperplane in $\A'$ before we meet $H_j$.
But we are moving parallel to every hyperplane in $\A_0$, so the basic hyperplane we must cross is $H_i$.
Thus $H_j$ follows $H_i$ in the ordering on $\A_1$, or in other words, $i<j$.
Since moving along arrows in $\D(\A,B)$ always moves us further in the ordering on $\A_1$, we can in particular never close a cycle.
\end{proof}

By induction, when $\A$ is supersolvable and $B$ is a canonical base region, we can cover $\D(\A,B)$ with $k$ acyclic induced sub-digraphs,
where $k$ is the rank of~$\A$.
Since this is exactly the lower bound of Proposition~\ref{lower}, we have given a tidier proof of the following theorem which was first 
proven in~\cite{hyperplane}.
\begin{theorem}
\label{dimension}
The order dimension of the poset of regions (with respect to a canonical base region) of a supersolvable hyperplane arrangement is equal to the
rank of the arrangement.
\qed
\end{theorem}
The proof of Proposition~\ref{superindependent} shows that if we order $\A_1$ as $H_1,H_2,\ldots,H_k$, we can construct the map
$\eta_{\A_1}$ of Lemma~\ref{binary}.
By induction, we obtain an explicit embedding in connection with Theorem \ref{dimension}.
A map very similar to $\eta_{\A_1}$ was considered in~\cite{hyperplane}, but an explicit embedding was not given there because of the lack of 
Proposition~\ref{box}. 

It is also possible to give a zonotopal embedding of the poset of regions (with respect to a canonical base region) of a supersolvable hyperplane arrangement.
\begin{theorem}
\label{super zone}
Let~$\A$ be a supersolvable hyperplane arrangement of rank $d$, and let~$B$ be a canonical base region.
Then $\Po(\A,B)$ has a zonotopal embedding in $\reals^d$.
\end{theorem}
\begin{proof}
Think of~$\A$ as a sequence $\A_1\subset\A_2\subset\cdots\subset \A_d=\A$ of supersolvable arrangements with $\rank(\A_i)=i$ and such that for each 
$i\in[d-1]$, Theorem~\ref{characterization} gives the partition $\A_i=\A_{i-1}\uplus(\A_i-\A_{i-1})$.
Since the canonical base region~$B$ was chosen according to an inductive definition, we have a canonical base region $B_i$ for each $\A_i$.
Choose $b_i$ to be a vector in $\Span(\A_i)\cap(\cap\A_{i-1})$ and choose the direction of $b_i$ so that, starting in $B_i$ and traveling in the direction of $b_i$,
one would reach the other $\A_i$-regions contained in $B_{i-1}$.
The vectors $b_i$ are used to define the components of the map $Z$, as defined in Section~\ref{zonotopal}.
Choose the directions of the normal vectors to~$\A$ as in Section~\ref{zonotopal}. 

We will prove by induction on $d$ that the normal vectors can be scaled so that for every $i\in[d]$ and every $H\in\A_i$, we have
\begin{equation}
\label{goal}
(n_H)_i>\sum_{H'\in \nu(H)}(n_{H'})_i.
\end{equation}
Then in particular, by Proposition~\ref{sufficient}, the map $Z$ defined in section~\ref{zonotopal} is a zonotopal embedding of $\Po(\A,B)$.
The case $d=1$ is trivial, so suppose $d\ge 2$, and consider first the case $i=d$ and then the case $i<d$.
For every $H'\in \A_{d-1}$, we have $(n_{H'})_d=0$ because $b_d\in H'$.
By Proposition~\ref{superindependent}, $\A_d$ induces an acyclic digraph of $\D(\A,B)$, so we can satisfy Inequality (\ref{goal}) with $i=d$
for every $H\in\A_d$.
In the case $i<d$, by induction we have for each $H\in\A_{d-1}$,
\[(n_H)_i>\sum_{H'\in \nu(H)\cap\A_{d-1}}(n_{H'})_i.\]
To satisfy Inequality (\ref{goal}) for each $i<d$ and each $H$, we need to be able to add into the right sides some terms arising from hyperplanes 
in $\A_d$.
Since the inequality is strict, this can be done as long as all the new terms are small enough.
To this end, we uniformly scale the normals to hyperplanes, preserving their relative proportions, and thus preserving Inequality (\ref{goal}) in the case $i=d$ as well.
\end{proof}

\section{Comments and questions}

\subsection*{The exceptional types}
The most immediate problem left unsolved is to determine the order dimension of the groups $E_6$, $E_7$, $E_8$, $F_4$, and $H_4$.
Absent further theoretical advances, this promises to be a computationally intense problem.
If any of the dimensions exceeds the rank of the arrangement, it would be the first example known to the author of a simplicial arrangement in 
which the dimension of the poset of regions exceeds the rank.
If each dimension is equal to the rank, is there a uniform proof of that fact (i.e.\ not relying on the classification of finite Coxeter groups)?

\subsection*{Quotients}
As noted in the introduction, Flath~\cite{Flath} determined the order dimension of the weak order on type A.
More generally, she determined the weak order for arbitrary (one-sided) quotients (with respect to parabolic subgroups) of the weak order 
on type A.
What are the dimensions of the quotients in other types?

\subsection*{Computation}
To embed the poset of regions by the method of Theorem~\ref{acyclic}, one needs to know the separating set of each element.
However, Theorem~\ref{acyclic} does lead to an improvement in computation.
Suppose that one wishes answer the question ``Is $R_1\le R_2$ in $\Po(\A,B)$?''
Suppose also that the basic unit of computation is to compute the answers to the questions ``Is $H$ in $S(R_1)$?'' and ``Is $H$ in $S(R_2)$?'' 
for a single $H\in\A$.
If at any point in the computation we get the answers ``yes'' and ``no'' to the two questions, we can conclude that $R_1\not\le R_2$.
If we begin with a covering of $\D(\A,B)$ by acyclic sub-digraphs $I_1,\ldots,I_d$ and test the hyperplanes within each sub-digraph in the 
order specified by Lemma~\ref{binary}, we obtain a further reduction:
Whenever we get the answers ``no'' and ``yes'' for a hyperplane $H\in I_k$, we can conclude that $\eta_{I_k}(R_1)<\eta_{I_k}(R_2)$, and it is 
not necessary to test the remaining hyperplanes in $I_k$.
This computational savings derives from ordering the hyperplanes in $I_k$ in a way that is compatible with $\D(\A,B)$, and possibly there is 
a more general computational scheme which is directly based on $\D(\A,B)$ or some variant.

\section{Acknowledgments}
The author wishes to thank Vic Reiner for helpful conversations and John Stembridge for helpful conversations and for writing computer programs 
to handle the exceptional types, as well as an anonymous referee for pointing out an error in a previous version of the proof of 
Proposition~\ref{j sigma}.

\newcommand{\journalname}[1]{\textrm{#1}}
\newcommand{\booktitle}[1]{\textrm{#1}}


\begin{thebibliography}{9}

\bibitem{Bar-Ihr}
H. Barcelo and E. Ihrig,
\textit{Modular elements in the lattice $L(A)$ when $A$ is a real reflection arrangement}, 
Selected papers in honor of Adriano Garsia, 
\journalname{Discrete Math.} {\bf 193} (1998), no.~1-3, 61--68.

\bibitem{BEZ}
A. Bj\"orner, P. Edelman and G. Ziegler,
\textit{Hyperplane Arrangements with a Lattice of Regions},
\journalname{Discrete Comput. Geom.} {\bf 5} (1990), 263--288.

\bibitem{Bourbaki}
N. Bourbaki,
\'{E}l\'{e}ments de math\'{e}matique. 
Groupes et alg\`{e}bres de Lie. Chapitres 4, 5 et 6. 
Masson, Paris, 1981.

\bibitem{boundedref}
N. Caspard, C. Le Conte de Poly-Barbut and M. Morvan,
\textit{Cayley lattices of finite Coxeter groups are bounded},
preprint, 2001.

\bibitem{Du-Mil}
B. Dushnik and E. Miller,
\textit{Partially ordered sets}, 
\journalname{Amer. J. Math.} {\bf 63} (1941), 600--610. 

\bibitem{Edelman}
P. Edelman,
\textit{A Partial Order on the Regions of $\reals^n$ Dissected by Hyperplanes},
\journalname{Trans. Amer. Math. Soc.} {\bf 283} no.~2 (1984), 617--631.

\bibitem{Fels-Trot}
S. Felsner and W. Trotter,
\textit{Dimension, Graph and Hypergraph Coloring},
\journalname{Order} {\bf 17} (2000) no.~2, 167--177.

\bibitem{Flath}
S. Flath,
\textit{The order dimension of multinomial lattices},
\journalname{Order} {\bf 10} (1993), no.~3, 201--219. 

\bibitem{FHRT}
Z. F\"{u}redi, P. Hajnal, V. R\"{o}dl and W. Trotter,
\textit{Interval orders and shift graphs}, in
\booktitle{Sets, graphs and numbers (Budapest, 1991)} 297--313, 
\journalname{Colloq. Math. Soc. János Bolyai} {\bf 60} (1992).

\bibitem{Humphreys}
J. Humphreys,
\booktitle{Reflection Groups and Coxeter Groups}, Cambridge Studies in Advanced Mathematics, {\bf 29},
Cambridge Univ. Press, 1990.

\bibitem{Or-Ter}
P. Orlik and H. Terao,
\booktitle{Arrangements of hyperplanes},
Grundlehren der Mathematischen Wissenschaften {\bf 300}, 
Springer-Verlag, 1992.

\bibitem{Rab-Riv}
I. Rabinovitch and I. Rival,
\textit{The Rank of a Distributive Lattice},
\journalname{Discrete Math.} {\bf 25} (1979) no.~3, 275--279.

\bibitem{hyperplane}
N. Reading
\textit{Lattice and Order Properties of the Poset of Regions in a Hyperplane Arrangement},
\journalname{Algebra Universalis}, to appear.

\bibitem{Trotter}
W. Trotter,
\booktitle{Combinatorics and Partially Ordered Sets:  Dimension Theory},
Johns Hopkins Series in the Mathematical Sciences,
The Johns Hopkins Univ. Press, 1992.

\bibitem{Yan}
M. Yannakakis,
\textit{The complexity of the partial order dimension problem},
\journalname{SIAM J. Algebraic Discrete Methods} {\bf 3} (1982) no.~3, 351--358. 

\bibitem{Ziegler}
G. Ziegler,
\textit{Combinatorial construction of logarithmic differential forms},
\journalname{Adv. Math.} {\bf 76} (1989), no.~1, 116--154. 


\end{thebibliography}
\end{document}